\documentclass[12pt]{amsart}
\usepackage{amsmath}
\usepackage{amsthm
}\usepackage{amssymb}
\usepackage{amscd}

\newcommand{\gothic}{\mathfrak}

\newcommand{\m}{{\gothic{m}}}

\newcommand{\Spec}{\operatorname{Spec}}

\newcommand{\Ann}{\operatorname{Ann}}
\newcommand{\height}{\operatorname{ht}}
\newcommand{\sdim}{\operatorname{sdim}}

\renewcommand{\phi}{\varphi}
\renewcommand{\to}{{\longrightarrow}}

\newcommand{\soc}{\operatorname{soc}}
\newcommand{\limt}{\lim\limits_{\longrightarrow t}}

\newtheorem{thm}{Theorem}

\newtheorem{prop}[thm]{Proposition}
\newtheorem{lemma}[thm]{Lemma}
\newtheorem{defn}[thm]{Definition}

\begin{document}
\title{The $F$-signature and strong $F$-regularity}
\author{Ian M. Aberbach}
\address{Mathematics Department \\
 	University of Missouri\\
	Columbia, MO 65211 USA}
\email{aberbach@math.missouri.edu}
\urladdr{http://www.math.missouri.edu/people/iaberbach.html}

\author{Graham J. Leuschke}
\address{Department of Mathematics \\
 	University of Kansas \\
 	Lawrence, KS	66045  USA} 
\email{gleuschke@math.ukans.edu}
\urladdr{http://www.math.ukans.edu/\textasciitilde gleuschke}
\date{\today}
\thanks{Aberbach was partially supported by the National Science Foundation and the National Security Agency.  Leuschke was partially supported by an NSF Postdoctoral Fellowship.}

\bibliographystyle{amsplain}

\numberwithin{thm}{section}

\begin{abstract} We show that the $F$-signature of a local ring of characteristic $p$, defined by Huneke and Leuschke, is positive if and only if the ring is strongly $F$-regular.
\end{abstract}  

\maketitle

In \cite{Huneke-Leuschke:2002}, Huneke and Leuschke define the $F$-signature of an $F$-finite local ring of prime characteristic with perfect residue field.  The $F$-signature, denoted $s(R)$, is an asymptotic measure of the proportion of $R$-free direct summands in a direct-sum decomposition of $R^{1/p^e}$, the ring of $p^e$th roots of $R$.  This proportion seems to give subtle information on the nature of the singularity defining $R$.  For example, the $F$-signature of any of the two-dimensional quotient singularities ($A_n$), ($D_n$), ($E_6$), ($E_7$), ($E_8$) is the reciprocal of the order of the group $G$ defining the singularity \cite[Example 18]{Huneke-Leuschke:2002}.  The main theorem of \cite{Huneke-Leuschke:2002} on $F$-signatures is as follows.

\begin{thm}\cite[Theorem 11]{Huneke-Leuschke:2002} Let $(R, \m, k)$ be a reduced complete $F$-finite Cohen--Macaulay local ring containing a field of prime characteristic $p$.  Assume that $k$ is perfect.  Then
\begin{enumerate}
\item If $s(R) >0$, then $R$ is weakly $F$-regular.
\item If in addition $R$ is Gorenstein, then $s(R)$ exists, and is positive if and only if $R$ is weakly $F$-regular.
\end{enumerate}\end{thm}

(See below for definitions of the $F$-signature and weak $F$-regularity.)  In this note, we extend this theorem in two directions: we remove the assumption in (2) that $R$ be Gorenstein, and we replace ``weakly $F$-regular" by ``strongly $F$-regular" throughout.  Our main theorem is thus as follows.

\begin{thm}\label{mainthm} Let $(R, \m,k)$ be a reduced excellent $F$-finite  local ring containing a field of characteristic $p$, and let $d=\dim R$.  Then the following are equivalent:
\begin{enumerate}
\item $\liminf \frac{a_q}{q^{d+\alpha(R)}} > 0$.
\item $\limsup \frac{a_q}{q^{d+\alpha(R)}} > 0$.
\item $R$ is strongly $F$-regular.
\end{enumerate}
 In particular, if the $F$-signature $s(R)$ is known to exist, then $s(R)$ is positive if and only if $R$ is strongly $F$-regular.\end{thm}

We also extend the definition of the $F$-signature to the case of an imperfect residue field.  This allows us to prove that $s(R)$ behaves well with respect to localization (Proposition~\ref{localization}).

Our results do not address the existence of the limit defining $s(R)$.  Yao has shown that $s(R)$ exists whenever $R$ is Gorenstein on the punctured
spectrum \cite{Yao:2002}.  

\section{The Main Result}

Throughout what follows, $(R, \m, k)$ is a reduced Noetherian local ring of dimension $d$, containing a field of positive characteristic $p$.  We use $q$ to denote a varying power of $p$. Set $d=\dim(R)$ and $\alpha(R) = \log_p[k:k^p]$.  We assume throughout that $R$ is $F$-finite, that is, the Frobenius endomorphism $F: R\to R$ defined by $F(r)=r^p$ is a module-finite ring homomorphism.  Equivalently, for each $q=p^e$, $R^{1/q} = \{r^{1/q}\ |\ r \in R\}$ is a finitely generated $R$-module.  In particular, this implies that $\alpha(R) < \infty$, and that $R$ is excellent \cite[Propositions 1.1 and 2.5]{Kunz:1976}.   Also, when computing length over $R$, we have $\lambda(R/I^{[q]}R) = \lambda(R^{1/q}/I R^{1/q})/q^{\alpha(R)}$.

We first define the $F$-signature of $R$. 

\begin{defn} Let $(R, \m ,k)$ be as above.  For each $q=p^e$, decompose $R^{1/q}$ as a direct sum of finitely generated $R$-modules $R^{a_q} \oplus M_q$, where $M_q$ has no nonzero free direct summands.  The {\it $F$-signature} of $R$ is
$$s(R) = \lim_{q\to \infty} \frac{a_q}{q^{d+\alpha(R)}},$$
provided the limit exists.\end{defn}

Our formulation differs slightly from the original definition in \cite{Huneke-Leuschke:2002}, where it is assumed that $k$ is perfect, or equivalently that $\alpha(R)=0$.  This reformulation allows us to show that $s(R)$ cannot decrease upon localization.  We use a lemma due to Kunz (\cite{Kunz:1976}).

\begin{lemma}\label{Kunz-alpha} Let $R$ be an $F$-finite Noetherian ring of characteristic $p$.  Then for any prime ideals $P \subseteq Q$ of $R$, $[k(P):k(P)^p]=[k(Q):k(Q)^p]p^{\dim R_Q/PR_Q}$.  In other words, $\alpha(R_P)=\alpha(R_Q)+\height{Q/P}$.  \end{lemma}

\begin{prop}\label{localization} Let $(R,\m)$ be an $F$-finite local ring and $P$ a prime ideal.  For $q=p^e$, let $a_q$ be the number of nonzero $R$-free direct summands in $R^{1/q}$, and let $b_q$ be the corresponding quantity for $R_P$.  Then
$$\frac{b_q}{q^{\dim(R_P) +\alpha(R_P)}} \geq \frac{a_q}{q^{\dim(R) +\alpha(R)}}.$$
In particular, if both $s(R)$ and $s(R_P)$ exist, then $s(R_P) \geq s(R)$.  
\end{prop}

\begin{proof} We have $(R_P)^{1/q} \cong (R^{1/q})_P$, so the number of $R_P$-free direct summands in $(R_P)^{1/q}$ is at least the number of $R$-free summands in $R^{1/q}$.  A straightforward computation using Lemma~\ref{Kunz-alpha} now gives the result.\end{proof} 

\medskip

We now begin to work toward showing that $s(R)$ is positive if and only if $R$ is strongly $F$-regular.  We refer the reader to \cite{Huneke:CBMS} for basic notions concerning the theory of tight closure, including {\it finitistic tight closure}, but review briefly the ideas used in the proof.  

A Noetherian ring $R$ of characteristic $p$ is said to be {\it weakly $F$-regular} provided every ideal of $R$ is tightly closed.  Equivalently, the zero module is finitistically tightly closed in $E=E_R(k)$, the injective hull of the residue field of $R$. In other symbols, $0_E^{*fg}=0$.  We say that $R$ is {\it strongly $F$-regular} if for every $c \in R$ not in any minimal prime of $R$, the inclusion $Rc^{1/q} \subset R^{1/q}$ splits for $q \gg 0$.  Equivalently, the zero module is tightly closed in $E$, that is, $0_E^* = 0$.  Weak and strong $F$-regularity are conjecturally equivalent, but this is known only in low dimension and in some special cases.

A {\it test element} for $R$ is an element $c$, not in any minimal prime of $R$, such that $cI^*\subseteq I$ for every ideal $I$ of $R$, and the {\it test ideal}, denoted $\tau(R)$,  is the ideal generated by all test elements.  For a reduced local ring $R$, $\tau(R) = \Ann_R 0^{*fg}_E$ by \cite[Theorem 8.23]{Hochster-Huneke:1990}.  Thus $R$ is weakly $F$-regular if and only if $\tau(R) = R$.  On the other hand, the {\it CS test ideal}, {\it cf.} \cite{Lyubeznik-Smith:2001} and \cite{Aberbach-Enescu:2002}, is the ideal $\tilde{\tau}(R) = \Ann_R 0^*_E$.  By work of \cite{Lyubeznik-Smith:2001} and \cite{Aberbach-Enescu:2002}, the CS test ideal behaves well under localization, so defines the non-strongly $F$-regular locus of $\Spec(R)$.  In particular, $R$ is strongly $F$-regular if and only if $\tilde{\tau}(R) = R$.

It is known that a weakly $F$-regular ring is {\it $F$-pure}, that is, the Frobenius morphism is  a pure homomorphism, and that for an $F$-pure ring both $\tau(R)$ and $\tilde{\tau}(R)$ are radical ideals.

A local ring $(R,\m, k)$ is said to be {\it approximately Gorenstein} provided there is a sequence $\{I_t\}$ of $\m$-primary irreducible ideals cofinal with the powers of $\m$.  When $R$ is Cohen--Macaulay and has a canonical ideal $J$ (so is Gorenstein at all associated primes), such a family can be obtained as follows:  Let $x_1, \ldots, x_d$ be a system of parameters such that $x_1 \in J$ and $x_2, \ldots, x_d$ form a system of parameters for $R/J$.  Then $I_t := (x_1^{t-1}J, x_2^t, \ldots, x_d^t)R$, for $t\geq 1$, gives the required family.  Furthermore, the direct limit $\limt R/I_t$, where the maps in the direct system are $R/I_t \overset{x_1\cdots x_d}\longrightarrow R/I_{t+1}$, is isomorphic to $E_R(k)$.  If $u_1 \in R$ is a representative for the socle generator of $R/I_1$, then $u_t := (x_1\cdots x_d)^{t-1}u_1$ generates the socle of $R/I_t$, and each $u_t$ maps in the limit to $u$, the socle element of $E_R(k)$.

More generally (\cite[Thm. 1.7]{Hochster:1977}), if $R$ is any  locally excellent Noetherian ring that is locally Gorenstein at associated primes, then $R$ is approximately Gorenstein.

The following result of Hochster, together with its corollary below, explains our interest in approximately Gorenstein rings. It can be thought of as  a generalization of \cite[Lemma 12]{Huneke-Leuschke:2002}.

\begin{prop}\label{H-splitting-criterion}\cite[Theorem 2.6]{Hochster:1977} Let $(R, \m)$ be an approximately Gorenstein local ring and let $\{I_t\}$ be a sequence of irreducible ideals cofinal with the powers of $\m$.  Let $f: R \to M$ be a homomorphism of finitely generated $R$-modules.  Then $f$ is a split injection if and only if $f \otimes_R R/I_t$ is injective for every $t$.
\end{prop}

\begin{prop}\label{splitting-criterion}  Let $(R, \m)$ be an approximately Gorenstein local ring with a family of irreducible ideals $\{I_t\}$ as above, and let $u_t \in R$ represent a socle generator for $R/I_t$.  Let $f:R \to M$ be a homomorphism of finitely generated $R$-modules. If $M$ has no free summands, then there exists $t_0>0$ such that $u_tM \subseteq I_tM$ for all $t \geq t_0$.
\end{prop}

\begin{proof} By Proposition~\ref{H-splitting-criterion}, $f \otimes R/I_t$ fails to be injective for some $t$.  Since $u_t$ is the unique socle element of $R/I_t$, we have $f(u_t) \in I_tM$, that is, $u_tM \subseteq I_tM$.  This continues to hold for all $t' \geq t$, since there is an injection $R/I_t \to R/I_{t'}$ with $u_t \mapsto u_{t'}$.
\end{proof}

We also use a result of Aberbach, which says that, in some sense, elements not in tight closures are very far from being in Frobenius powers.

\begin{thm}\cite[Prop. 2.4]{Aberbach:2001}\label{big-powers} Let $(R,\m)$ be an excellent local domain such that the completion is also a domain.  Let $N = \limt R/J_t$ be a direct limit system of cyclic modules.  Fix $u \not\in 0^*_N$.  Then there exists $q_0$ such that 
$$\bigcup_t(J^{[q]}_t:u^q_t) \subseteq \m^{[q/q_0]}$$
 for all $q \gg 0$ (where the sequence $\{u_t\}$ represents $u\in N$ and $u_t \mapsto u_{t+1}$).   
\end{thm}

\medskip

\begin{proof}[Proof of Theorem~\ref{mainthm}] 
The Cohen-Macaulayness of $R$ is forced by the assumptions (\cite[Theorem
11]{Huneke-Leuschke:2002} and \cite{Hochster-Huneke:1990}), so we may assume
throughout that $R$ is Cohen-Macaulay.

That (1) implies (2) is trivial.  So assume that (2) holds.
 We proceed by induction on the dimension $d$, the case $d=0$ being trivial.  If $d>0$, then Proposition~\ref{localization} shows that we may assume by induction on $d$ that $R$ is strongly $F$-regular on the punctured spectrum.  We will show that $0^*_E = 0$, where as above $E=E_R(k)$ is the injective hull of the residue field of $R$.  

Since $\tilde{\tau}(R)=\Ann_R 0^*_E$ is a radical ideal and is known to define the non-strongly $F$-regular locus of $R$ (see \cite{Aberbach-Enescu:2002}), and $R$ is strongly $F$-regular on the punctured spectrum, $\Ann_R 0^*_E$ contains the maximal ideal $\m$. If $\tilde{\tau}(R)=R$, then we are done, so we assume $\tilde{\tau}(R)=\m$.  Then $0^*_E = \soc(E)$.

As in the discussion above, $E=E_R(k) \cong \limt R/I_t$ for a family of irreducible ideals $I_t$.  Let $u$ be a socle generator for $E$ and $\{u_t\} \subseteq R$ a sequence of representatives for the socle generators of $R/I_t$, converging to $u$.

Fix a power $q$ of the characteristic, and decompose $R^{1/q} \cong R^{a_q} \oplus M_q$, where $M_q$ has no nonzero free summands. Then for each $t$, we have
\begin{equation*}
\begin{split}
\lambda\bigg(R/I_t^{[q]}\bigg) - \lambda\bigg(R/(I_t,u_t)^{[q]}\bigg)
        &=\frac{\lambda\left(R^{1/q}/I_t R^{1/q}\right)}{q^{\alpha(R)}} - \frac{\lambda\left(R^{1/q}/(I_t,u_t)R^{1/q}\right)}{q^{\alpha(R)}} \\
        &= \frac{a_{q}\lambda(R/I_t) + \lambda(M_q/I_t M_q)}{q^{\alpha(R)}} \\
        & \phantom{blergh} -  \frac{a_{q} \lambda\left(R/(I_t,u_t)\right) -
         \lambda\left(M_q/(I_t,u_t)M_q\right)}{q^{\alpha(R)}}  \\
        &= \frac{a_{q}\lambda(R/I_t) - a_q \lambda(R/(I_t, u_t))}{q^{\alpha(R)}}\\
        & \phantom{blergh} + \frac{\lambda(M_q/I_tM_q) - \lambda(M_q/(I_t,u_t)M_q)}{q^{\alpha(R)}} \\
        &= \frac{a_{q}+c_{t,q}}{q^{\alpha(R)}},
\end{split} 
\end{equation*}
for some $c_{t,q} \geq 0$.  By Proposition~\ref{splitting-criterion}, there exists $t_0 > 0$ such that $u_{t} M_q \subseteq I_{t} M_q$ for $t\geq t_0$, that is, $c_{t,q}=0$ for $t\geq t_0$.  On the other hand, 
$\lambda(R/I_t^{[q]}) - \lambda(R/(I_t,u_t)^{[q]}) = \lambda(R/(I_t^{[q]}:u_t^{q}))$ is equal to $1$ for large $t$ since $(I_t^{[q]}:u_t^q) = \m$ for large $t$.
Thus, for large $t$,
$$\lim_{q\to \infty} \frac{a_q+ c_{t,q}}{q^{d+\alpha(R)}} = \lim_{q\to\infty}\frac{1}{q^{d+\alpha(R)}} = 0,$$
a contradiction.

Lastly, assume that $R$ is strongly $F$-regular and keep the same notation.  We then have $0^*_E=0$, so $u\not\in 0^*_E$. By Theorem~\ref{big-powers}, then, there exists $q_0$ such that 
$$(I_t^{[q]} :_R u_t^q) \subseteq 
\m^{[q/q_0]}$$
for all $q \geq q_0$.  
Fix $q \geq q_0$.  Then there exists $t_0$ such that for all $t \geq t_0$ we have 
\begin{equation*}
\begin{split}\frac{a_q}{q^{\alpha(R)}}
	&= \lambda\left(R/I_t^{[q]}\right) - \lambda\left(R/(I_t^{[q]}, u_t^q)\right)\\
	&= \lambda\left(I_t^{[q]}: u_t^q\right) \\
	&\geq \lambda\left(R/\m^{[q/q_0]}\right).
\end{split}
\end{equation*}
Divide by $q^d$ and pass to the limit; we see that $\liminf \frac{a_q}{q^{d+\alpha(R)}}
 \geq \operatorname{e}_{HK}(\m,R)/q_0^d > 0.$   Thus (1) holds.

The last statement is immediate if there is a limit.
\end{proof}

\medskip

The $F$-signature suggests a form of dimension that we may attach to
an $F$-finite reduced local ring.  Let $s_j = \lim_{q\to\infty} \dfrac{a_q}{q^{j+\alpha(R)}}$ for $0 \le j \le d = \dim(R)$ and set  $s_{-1} =1$. Then we can define the $s$-dimension of $R$ as $\sdim(R) = \max \{ j \ge -1| s_j > 0\}$.  A ring which is $F$-pure then has non-negative $s$-dimension,
and  Theorem~\ref{mainthm} says that $R$ is strongly $F$-regular if and
only if $\sdim(R) = \dim(R)$.


\providecommand{\bysame}{\leavevmode\hbox to3em{\hrulefill}\thinspace}

\end{document}